\documentclass[11pt]{article} 
\usepackage{amssymb}
\usepackage{latexsym,amsmath,amsthm,amsfonts}

\newtheorem{prop}{Proposition}

\newtheorem{thm}{Theorem}
\newtheorem{cor}{Corollary}

\newtheorem{lemref}{Lemma}

\newcommand{\be}{\begin{equation}}
\newcommand{\ee}{\end{equation}}

\def\real{\mathbb R}

\def\P{\mathbb P}

\newcommand{\B}{\mbox{$\cal B$}}
\newcommand{\C}{\mathcal{C}}
\newcommand{\sigS}{\mbox{$\cal S$}}

\newcommand{\F}{\mbox{$\mathcal F$}}
\newcommand{\G}{\mbox{$\cal G$}}

\newcommand{\X}{\mbox{$\cal X$}}

\newcommand{\bX}{\mbox{$\bf X$}}

\begin{document}

\bibliographystyle{plain}

\title{Uniform Approximation of Vapnik-Chervonenkis Classes}

\author{Terrence M.\ Adams
\thanks{Terrence Adams is with the Department of Defense, 
9800 Savage Rd. Suite 6513, Ft. Meade, MD 20755} \ 
and Andrew B.\ Nobel 
\thanks{Andrew Nobel is with the Department of Statistics and 
Operations Research, University of North Carolina, Chapel Hill,
NC 27599-3260.  Email: nobel@email.unc.edu}}

\date{September 2010}

\maketitle

\begin{abstract}
For any family of measurable sets in a probability
space, we show that either (i) the family has infinite Vapnik-Chervonenkis (VC) 
dimension or (ii) for
every $\epsilon > 0$ there is a finite partition $\pi$
such the $\pi$-boundary of each set has measure 
at most $\epsilon$.
Immediate corollaries include 
the fact that a family with finite VC dimension has finite bracketing 
numbers, and satisfies
uniform laws of large numbers for every ergodic process.  
From these corollaries, we derive analogous results for
VC major and VC graph families of functions.  
\end{abstract}

\newpage

\section{Introduction}


Let $(\X,\sigS,\mu)$ be a probability space and let $\C \subseteq \sigS$ 
be a given family of measurable sets.  The Vapnik-Chervonenkis dimension
of $\C$ is a measure of its combinatorial complexity, specifically, the ability
of $\C$ to separate finite sets of points.   
Given a finite set $D \subseteq \X$,
let $\{ C \cap D : C \in \C \}$ be the collection of subsets of
$D$ selected by the members of $\C$.  The family $\C$ is
said to shatter $D$ if its elements can select every subset of $D$, or
equivalently, if $| \{ C \cap D : C \in \C \} | =  2^{|D|}$.
Here and in what follows, $| A |$ denotes the cardinality
of a given set $A$.  The Vapnik-Chervonenkis (VC) dimension \cite{VC71} of
$\C$, denoted $\dim(\C)$, is the largest integer $k$ such that 
$\C$ is able to shatter {\em some} set of cardinality $k$.  
If $\C$ can shatter arbitrarily large finite sets, then $\dim(\C) = +\infty$.  
A family of sets $\C$ is said to be a VC class if $\dim(\C)$ is finite.

Let $\pi$ be a finite, measurable partition of $\X$.  For every set $C \in \C$,
the $\pi$-boundary of $C$, denoted $\partial(C : \pi)$, is the union of all the
cells in $\pi$ that intersect both $C$ and its complement
with positive probability.  Formally,
\[
\partial(C : \pi)
\ = \ 
\cup \, \{ A \in \pi : \mu(A \cap C) \cdot  \mu(A \cap C^c) > 0 \} .
\]
Note that $\partial(C : \pi)$ depends on $\mu$, though this dependence
is suppressed in our notation.  We will call a family $\C$ 
{\em finitely approximable}
if for every $\epsilon > 0$ there exists a finite, measurable partition $\pi$ of
$\X$ such that $\mu(\partial(C : \pi) ) \leq \epsilon$ for every $C \in \C$.
Our principal result is the following.

\vskip.2in

\begin{thm}
\label{BdyThm}
Let $(\X, \sigS,\mu)$ be a probability space and let
$\C \subseteq \sigS$ be any family of sets.  Then either 
(i) $\C$ is finitely approximable or (ii) $\C$ has infinite VC dimension.
\end{thm}

\vskip.2in

Theorem \ref{BdyThm} extends immediately to finite positive measures;
we restrict attention to the case of probability measures for simplicity.
Gaenssler and Stute \cite{GaenStute76} studied $\pi$-boundaries in work
on uniform convergence of measures.  In conjunction with Theorem
\ref{BdyThm}, their results show that, if for some VC-class $\C$ and 
some sequence $\{ \mu_n \}$
of finite measures, $\mu_n(A) \to \mu(A)$ for every $A \in \sigma(\C)$,
then this convergence is uniform over $\C$.  One may establish
the same conclusion using Corollary \ref{Bracket}.

In general, alternatives (i) and (ii) of Theorem \ref{BdyThm} are not
mutually exclusive: there exist families $\C$ that are finitely approximable
and have infinite VC dimension.  Moreover the finite approximability of 
$\C$ will generally depend on the measure $\mu$.  
To take a simple example, let $\C$ be
the family of all
Borel measurable subsets of the unit interval $[0,1]$.  Then $\C$ clearly
has infinite VC dimension.  An easy argument shows that $\C$ is finitely 
approximable if $\mu$ has countable support, but that $\C$ is not
finitely approximable if $\mu$ is absolutely continuous with respect to
Lebesgue measure.  As the following, equivalent, version of Theorem
\ref{BdyThm} makes clear, families with finite VC dimension are
finitely approximable for any probability measure $\mu$.

\vskip.1in

\begin{thm}
\label{BdyThm2}
Let $(\X, \sigS)$ be a measurable space. 
If $\C \subseteq \sigS$ has finite VC dimension, then $\C$ is finitely
approximable for any probability measure $\mu$.
\end{thm}

\vskip.1in

Families of sets with finite VC-dimension figure prominently in 
machine learning, empirical process theory and combinatorial
geometry ({\it c.f.} 
\cite{Poll84, VaaWel96, DevGyoLug96, Dudley99, Vapnik00, Mat02})
and have been widely studied in these fields.  The majority of this
work concerns the combinatorial properties of VC-classes, and related
exponential probability inequalities for uniform laws of large numbers
under independent sampling (see Section \ref{ULLN} below).  The uniform
approximation guaranteed by Theorem
\ref{BdyThm2} provides new insights into the structure of VC-classes.  

Some immediate corollaries of Theorem \ref{BdyThm2} are explored 
in Sections \ref{Brack} and \ref{ULLN} below, including new results
on the bracketing properties of VC major and VC graph classes of 
functions.  Approximation properties analogous to those of 
Theorems \ref{BdyThm} 
and \ref{BdyThm2} may be established for classes of functions
with finite fat-shattering (gap) dimension \cite{KeaSch94}
by extending the arguments in Section \ref{PfBdyThm}.  

\vskip.1in

The proof of Theorem \ref{BdyThm} makes use of an equivalent version
of the VC dimension that we now describe.  Recall that the join of $k$ sets 
$A_1,\ldots, A_k \subseteq \X$, 
denoted $J = \bigvee_{i=1}^k A_i$, is the finite partition of $\X$ consisting of all 
non-empty intersections 
$\tilde{A}_1 \cap \cdots \cap \tilde{A}_k$,
where $\tilde{A}_i \in \{ A_i, A_i^c \}$ for $i = 1,\ldots,k$.  Equivalently, $J$
consists of the non-empty atoms of the field generated by $A_1,\ldots, A_k$. 
The collection $A_1,\ldots, A_k \subseteq \X$ is said to be Boolean independent
if $J$ has (maximal) cardinality $2^k$.   The dual VC dimension, denoted
$\dim^*(\C)$, is the largest $k$ such that $\C$ contains $k$ Boolean independent
sets.  If $\C$ contains Boolean independent families of every finite size,
then $\dim^*(\C) = +\infty$.  The dual VC-dimension
was introduced by Assouad \cite{Assouad83}, and is so named because
$\dim^*(\C)$ is the VC-dimension of the dual family 
$\{ D_x : x \in \X \} \subseteq 2^{\cal C}$, where $D_x = \{ C \in \C : x \in C \}$. 
We will make use of the following, elementary result, whose proof
can be found in \cite{Assouad83}, see also \cite{Mat02, AdNob10a}. 

\vskip.1in

\begin{lemref}
\label{Join}
Let $\C$ be any collection of subsets of $\X$.  The VC-dimension
$\dim(\C)$ is finite if and only if the dual VC-dimension $\dim^*(\C)$
is finite.
\end{lemref}

In proving Theorem \ref{BdyThm} we begin with the assumption that $\C$
is not finitely approximable, and then deduce from this that 
$\dim^*(\C) = + \infty$.  Specifically, we show that 
for every $L \geq 1$ the family $\C$ contains a sub-family of
$L$ Boolean independent sets.  We note that Boolean independence
plays a related role in work of Rosenthal \cite{Rosen74}, who shows that
if a sequence of sets $\{ C_n : n \geq 1\}$ contains no pointwise 
convergent subsequence,
then there is an infinite subsequence $\C_0 = \{ C_{n_m} : m \geq 1 \}$ such
that each finite subfamily of $\C_0$ is Boolean independent.  

The construction of Boolean independent sets in Theorem 1 
proceeds in stages.  At each
stage a splitting set is produced by means of a weak limit, and is then
incorporated in the construction of the splitting sets at subsequent stages.
The resulting sequence of splitting sets is used to identify
Boolean independent collections of arbitrary finite size.
As noted by Ramon van Handel (private communication),
the proof of Theorem 1 has points of intersection with
the construction of a critical set for product measures in
Theorem 11-1-1 of Talagrand \cite{Tala84}, and with the 
notion of weakly dense sequences in \v{C}ech-complete spaces 
employed by Bourgain, Fremlin, and Talagrand \cite{BourFremTala78}.
Essential differences emerge from a number of 
factors, including our focus on finite approximation under a fixed 
(but arbitrary) distribution in the absence
of topological structure, as well as the recursive 
construction of splitting sets that is employed in the theorem.

\subsection{Overview}

The next two sections are devoted to corollaries of Theorem \ref{BdyThm} 
to families of sets and functions with bounded combinatorial complexity. 
In Section \ref{Brack} we establish that VC classes of sets have finite bracketing
numbers, and deduce similar results for VC major and VC graph families of
functions.  In Section \ref{ULLN} we show that VC classes satisfy uniform
laws of large numbers for every ergodic process.
The proof of Theorem \ref{BdyThm} is presented in 
Section \ref{PfBdyThm}.

\section{Bracketing of VC Classes of Sets and Functions}
\label{Brack}

Let $\F$ be a family of measurable functions $f: \X \to \real$.  
We recall some basic definitions from the theory of empirical processes.
A measurable function $F: \X \to [0,\infty)$ 
is said to be an envelope for $\F$ if $|f(x)| \leq F(x)$ 
for each $x \in \X$ and $f \in \F$.  
The family $\F$ is said to be separable if there is a countable
sub-family $\F_0 \subseteq \F$ such that each function $f \in \F$
is a pointwise limit of a sequence of functions in $\F_0$.
For each pair of measurable functions $g, h: \X \to \real$
with $g \leq h$, the bracket $[g,h]$ denotes the set of all
measurable functions $f$ such that $g \leq f \leq h$ pointwise on $\X$.  
In particular, $[g,h]$ is said to be an $\epsilon$-bracket 
if $\int (h - g) d\mu \leq \epsilon$.
For $\epsilon > 0$, the bracketing number $N_{[ \, ]}(\epsilon,\F,\mu)$ 
of $\F$ is the least number 
of $\epsilon$-brackets needed to cover $\F$.  In general, the 
functions defining the minimal brackets need not be elements of $\F$.

\subsection{VC Classes of Sets}

Let a measure $\mu$ and family $\C \subseteq \sigS$ be fixed. 
The notions of separability and bracketing may be applied to $\C$
if we regard its elements as indicator functions.  In this
case we may assume, without loss of generality, that the lower
and upper limits of each bracket are themselves indicator functions.

\vskip.1in

\begin{cor}
\label{Bracket}
If $\C$ is a separable VC-class, then $N_{[ \, ]}(\epsilon,\C,\mu)$ is finite 
for every $\epsilon > 0$.
\end{cor}

\noindent
{\bf Proof:} 
By routine arguments, we may assume that $\C$ is countable.
Fix $\epsilon > 0$.  
Let $\pi = \{A_1,\ldots,A_m\}$ be a finite measurable 
partition of $\X$ such that $\mu(\partial(C : \pi)) < \epsilon$ 
for every $C \in \C$, and assume
without loss of generality that each cell of $\pi$ has positive
$\mu$-measure.  
For each $C \in \C$, remove all points in $C$ from $A_j$ if $\mu(A_j \cap C) = 0$,
and remove all points in $C^c$ from $A_j$ if  
$\mu(A_j \cap C^c) = 0$.  Denote the resulting
set by $B_j$.  Clearly $B_j \subseteq A_j$ and 
$\mu(A_j \setminus B_j) = 0$ as $\C$ is countable.  The definition
of $B_j$ ensures that for each $C \in \C$ exactly one of the
following relations holds: $B_j \subseteq C$, 
$B_j \subseteq C^c$, or $\mu(B_j \cap C) \cdot \mu(B_j \cap C^c) > 0$.
Let $B_0 = \X \setminus \cup_{j=1}^m B_j$, and define the partition
$\pi' = \{ B_0, B_1, \ldots, B_m \}$.  Given $C \in \C$ let
$C_l = \cup \{ B \in \pi' : B \subseteq C \}$ and  
$C_u = \cup \{ B \in \pi' : B \cap C \neq \emptyset \}$.
A straightforward argument shows that 
$C_l \subseteq C \subseteq C_u$, and that
$\mu(C_u \setminus C_l) = 
\mu(\partial(C : \pi')) = \mu(\partial(C : \pi)) < \epsilon$.  
It follows that $\Theta = \{ [C_l,C_u] : C \in \C\}$ is a collection of 
$\epsilon$-brackets covering $\C$.   The cardinality of $\Theta$
is at most $2^{2 | \pi' |}$.

\subsection{VC Major Families}

Let $\F$ be a family of measurable functions $f: \X \to \real$ with envelope $F$.
For $f \in \F$ and $\alpha \in \real$ let $L_f(\alpha) = \{ x : f(x) \leq \alpha \}$ 
be the $\alpha$-level set of $f$.  Define
\[
\C_\alpha = \left\{ \, L_f(\alpha) : f \in \F \, \right\}
\]
to be the family of $\alpha$-level sets associated with functions in $\F$.

\begin{prop}
\label{VCmajor}
Suppose that $\dim( \C_\alpha ) < \infty$ for every $\alpha \in \real$. 
If $\mu$ is any probability measure on $(\X, \sigS)$
such that $\int F \, d\mu < \infty$, then
$N_{[\,]}(\epsilon, \F, \mu) < \infty$ for every $\epsilon > 0$.
\end{prop}

\noindent
{\bf Proof:}
Suppose first that $\F$ is bounded, with constant envelope $M < \infty$.  
Fix $\epsilon > 0$ and let $K$ be an integer such that 
$2 M / K \leq \epsilon$.  For each $f \in \F$ define the approximation
\[
\tilde{f}(x)
\ = \ 
M - \frac{2M}{K} \sum_{j=1}^K I( x \in L_f(\alpha_j) )
\ \mbox{ with } \ 
\alpha_j = M - \frac{2M j}{K} .
\]
The choice of $M$ and $K$ ensure that
$\tilde{f}(x) - \epsilon \leq f(x) \leq \tilde{f}(x)$ for each $x \in \X$.  
The dimension of $\C_{\alpha_j}$ is finite by assumption,
and it then follows from Corollary \ref{Bracket} that there is a 
finite collection $\Theta_j$ of $\epsilon/2M$-brackets that 
covers the level sets $\{ L_f(\alpha_j) : f \in \F \}$.    
For each $f \in \F$ let $[g_f^j, h_f^j]$ be 
a bracket in $\Theta_j$ containing $L_f(\alpha_j)$.  With this
identification, define upper
and lower approximations of $f$ as follows:
\[
\tilde{f}_l 
\ = \ 
M - \frac{2M}{K} \sum_{j=1}^K h_f^j(x) - \epsilon 
\ \mbox{ and } \ 
\tilde{f}_u 
\ = \ 
M - \frac{2M}{K} \sum_{j=1}^K g_f^j(x) 
\]
An easy argument shows that $\tilde{f}_l \leq f \leq \tilde{f}_u$, and
the family of brackets $\Theta = \{ [ \tilde{f}_l,  \tilde{f}_u ] : f \in \F \}$
is finite, as $|\Theta| \leq \Pi_{j = 1}^K |\Theta_j|$.  Moreover,
\[
\tilde{f}_u - \tilde{f}_l 
\ \leq \ 
\frac{2M}{K} \sum_{j=1}^K (h_f^j(x) - g_f^j(x)) + \epsilon ,
\]
and therefore $\int(\tilde{f}_u - \tilde{f}_l) d\mu \leq 2 \epsilon$. 
Thus $\Theta$ is a finite family of $2 \epsilon$-brackets covering $\F$.

Suppose now that $\F$ has an envelope $F$ such that $\int F d\mu < \infty$.  
Given $\epsilon > 0$ let $M < \infty$ be such that 
$\int_{F > M} F d\mu < \epsilon$.
For each $f \in \F$ define the truncation $f_M(x) = (f(x) \vee -M) \wedge M$,
and let $\F_M = \{ f_M : f \in \F \}$.  By the preceding argument,
there is a finite family $\Theta$ of $\epsilon$-brackets covering 
$\F_M$.  Let $[g,h]$ be an element of $\Theta$; without loss of 
generality, we may assume that $|g|, |h| \leq M$.  Define
\[
g' \ = \ g \wedge (-F I(F > M)) \ \mbox{ and } \ 
h' \ = \ h \vee (F I(F > M) )
\]
and note that $g' \leq g \leq h \leq h'$.  Moreover, $f_M \in [g,h]$ implies
$f \in [g',h']$, so the finite family of brackets $\{ [g',h'] : [g,h] \in \Theta \}$ 
covers $\F$.   It is easy to see that 
\[
h' - g' = (h-g) I(F \leq M) + 2 FI(F > M) ,
\]
and therefore
$
\int(h' - g') d\mu 
\ \leq \ 
\int(h - g) d\mu + 2 \int_{F > M} F d\mu 
\leq 
3 \epsilon.
$

\subsection{VC Graph Families} 

Let $\F$ be a family of measurable functions $f: \X \to \real$ with envelope
$F(x)$.  The graph of $f \in \F$ is defined by 
\[
G_f = \{ (x,s) : x \in \X \, \mbox{ and } \, 0 \leq s \leq f(x) 
\, \mbox{ or } \, f(x) \leq s \leq 0 \} \,
\subseteq \X \times \, \real.
\]
Let ${\cal G}(\F) = \{ G_f : f \in \F \}$ be the family of graphs of
functions in $\F$.  

\vskip.4in

\begin{prop}
\label{VCgraph}
Suppose that $\dim({\cal G}(\F)) < \infty$.
If $\mu$ is any probability measure on $(\X, \sigS)$
such that $\int F \, d\mu < \infty$, then 
$N_{[\,]}(\epsilon, \F, \mu) < \infty$ for each $\epsilon > 0$.
\end{prop}

\noindent
{\bf Proof:}
Suppose first that $\F$ is bounded, with constant envelope $M < \infty$.  
The finiteness of the bracketing numbers is not affected if we
replace each function $f \in \F$ by $(f + M)/2M$, and we therefore
assume that every $f \in \F$ takes values in $[0,1]$.  With this restriction,
\[
G_f \ = \ \{ (x,s) : x \in \X \mbox{ and } 0 \leq s \leq f(x) \leq 1 \} 
\, \subseteq \, \X \times [0,1] .
\]
Let $\lambda(\cdot)$ denote Lebesgue measure on the Borel subsets
$\B$ of $[0,1]$, and
define the product measure $\nu = \mu \otimes \lambda$ on
$(\X \times [0,1], \sigS \otimes \B)$.

Fix $\epsilon > 0$.  As $\G(\F)$ has finite VC dimension, 
Corollary \ref{Bracket} ensures that $\G(\F)$ is covered
by a finite collection $\Theta$ of $\epsilon$-brackets.  
Without loss of generality, we may represent the brackets in 
$\Theta$ in the form $[A,B]$, where $A, B \in \sigS \otimes \B$
and $A \subseteq B$.
Let $[A,B]$ be a bracket in $\Theta$. 
For each $x \in \X$ define 
\[
g(x) = \mbox{ess-sup} (\{ s : (x,s) \in A \})
\ \mbox{ and } \ 
h(x) = \mbox{ess-inf} (\{ s : (x,s) \in B^c \}) ,
\]
where for $U \subseteq [0,1]$ the essential supremum
$\mbox{ess-sup}(U) = \inf\{ \alpha : \mu(U \cap [0,\alpha]) = \mu(U) \}$,
and $\mbox{ess-inf}(U)$ is defined analogously.  Routine arguments shows that
$g$ and $h$ are measurable, that $g \leq h$, and that
$\nu(A \setminus G_g) = \nu(B^c \setminus G_h^c) = 0$.  Moreover,
for every function $f:\X \to [0,1]$ it is easy to see that $G_f \in [A,B]$ implies  
$G_f \in [G_g,G_h]$, which implies in turn that $g \leq f \leq h$.

It follows from the arguments above that the finite family $\Theta_0$ 
of brackets $[g,h]$ derived from the elements of $\Theta$ covers $\F$.  
In order to assess the size of these brackets, note that 
\[
(G_h \setminus G_g)_x 
\ = \ 
\{ s : (x,s) \in G_h \setminus G_g \}
\ = \ 
\{ s : g(x) < s \leq h(x) \} 
\]
and therefore by
Fubini's theorem
\[
\int ( h(x) - g(x) ) d \mu(x)  
\ = \
\int \lambda( (G_h \setminus G_g)_x ) d \mu(x) 
\ = \ 
\nu( G_h \setminus G_g) 
\ \leq \ 
\nu(B \setminus A) 
\ \leq \ 
\epsilon .
\]
Thus every element $[g,h]$ of $\Theta_0$ is an $\epsilon$-bracket 
under $\mu$. 

The argument for an unbounded family $\F$ with an integrable 
envelope $F$ is similar to that for VC Major families.
Given $\epsilon > 0$ let $M < \infty$ be such that 
$\int_{F > M} F d\mu < \epsilon$.
For each $f \in \F$ define the truncation $f_M(x) = (f(x) \vee -M) \wedge M$,
and let $\F_M = \{ f_M : f \in \F \}$.  As
$G_{f_M} = G_f \cap (\X \times [-M,M])$, it is easy to see that the 
dimension of $\G(\F_M)$ is no greater than that of $\G(\F)$, and is
therefore finite.  The preceding argument shows that 
there is a finite collection of $\epsilon$-brackets covering $\F_M$,
and these can be extended to $3 \epsilon$-brackets covering $\F$
following the proof of Proposition \ref{VCmajor}.

\section{Uniform Laws of Large Numbers}
\label{ULLN}

Let $\bX = X_1, X_2, \ldots$ be a stationary ergodic 
process taking values in $(\X,\sigS)$.  The ergodic theorem 
ensures that for every measurable set $C$ the sample averages 
$n^{-1} \sum_{i=1}^n I_C(X_i)$ converge almost surely 
to $P(X \in C)$.  A family $\C \subseteq \sigS$ satisfies 
a uniform laws of large numbers with respect to $\bX$ if
the discrepancy
\[
\Delta_n(\C : \bX) \ = \ 
\sup_{C \in {\cal C}} 
\left| \frac{1}{n} \sum_{i=1}^n I_C(X_i) - P(X \in C) \right|
\]
tends to zero almost surely as $n$ tends to infinity, so that 
the relative frequencies 
of sets in $\C$ converge uniformly to their limiting probabilities.


For i.i.d.\ processes $\bX$, Vapnik and Chervonenkis \cite{VC71} 
gave necessary and sufficient conditions under which 
$\Delta_n(\C : \bX) \to 0$.  For VC-classes they
established exponential inequalities of the form 
$\P( \Delta_n(\C : \bX) > t ) \leq a \cdot n^{\dim(\cal C)} \cdot \exp\{-b t^2\}$,
where $a,b$ are positive constants independent of $\bX$ and $\C$.
Consequently, VC classes have uniform laws of large numbers for
any i.i.d.\ process.
Talagrand \cite{Tala87} provided necessary and sufficient conditions 
for uniform laws of large numbers that strengthen those of \cite{VC71}:
for non-atomic distributions, 
$\Delta_n(\C : \bX) \not\to 0$ if and only 
if there is a set $A \in \sigS$ with $P(A) > 0$ such that with probability
one $\C$ shatters every finite subset of $\{ X_i: X_i  \in A \}$.

Using the bracketing properties of VC classes established in the
previous section one may immediately extend this result
to the general ergodic case.
The following theorem appears in Adams and Nobel \cite{AdNob10a}
(under an additional Polish assumption), 
where there is also a discussion of related work on uniform
laws of large numbers under dependent sampling.

\begin{thm}
\label{Ergodic}
If $\C$ is a separable VC-class of sets and $\bX$ is a stationary ergodic 
process, then $\Delta_n(\C : \bX) \to 0$ almost
surely
as $n$ tends to infinity.
\end{thm}

\noindent
{\bf Proof:} The stated convergence follows easily from Corollary \ref{Bracket}
and standard arguments for the Blum DeHardt law of large numbers 
({\it c.f.} \cite{VaaWel96,Dudley99}).

\vskip.2in

One may establish uniform laws of large numbers for
separable VC major and VC graph classes of functions 
in the general ergodic case using the bracketing results
in Propositions \ref{VCmajor} and \ref{VCgraph}, respectively.  
In \cite{AdNob10a} these results are derived directly from 
Theorem \ref{Ergodic}. Related work for families of functions, 
under a more general, scale specific, notion of dimension can
be found in \cite{AdNob10b}.

\section{Proof of the Main Theorem}
\label{PfBdyThm}

In the case where $\X$ is a complete separable 
metric space and $\sigS$ is the Borel subsets
of $\X$, one may prove Theorem \ref{BdyThm} using
arguments similar to those used in \cite{AdNob10a} to establish uniform 
laws of large numbers for VC classes under ergodic 
sampling.  The details can be found in an earlier version 
\cite{AdNob10c} of the results presented here.  Below we provide
a simpler argument that does not require the Polish assumption.
The new argument, which follows the outline of the proof in
\cite{AdNob10a}, employs several simplifications and improvements that
were suggested by an anonymous referee of \cite{AdNob10a}, in particular, the
use of Hilbert space weak limits in the definition of splitting sets.

\subsection{Proof of Theorem \ref{BdyThm}} 

It follows from standard results on the $L_p$-covering
numbers of VC classes (for example, Theorem 2.6.4 of 
\cite{VaaWel96}) that there exists a countable sub-family
$\C_0$ of $\C$ such that
$
\inf_{C' \in {\cal C}_0} \mu(C' \triangle C) = 0
$
for each $C \in \C$.
An elementary argument then shows that
\[
\sup_{C \in {\cal C}} \mu( \partial(C : \pi)) \ = \ 
\sup_{C \in {\cal C}_0} \mu( \partial(C : \pi)) 
\]
for every finite partition $\pi$,
and we may therefore assume that $\C$ is countable.  
Let $\C = \{C_1, C_2, \ldots \}$ and let 
$\sigS_0 = \sigma(\C) \subseteq \sigS$ be the sigma field generated
by $\C$. 
Suppose that the uniform approximation property fails to hold for 
$\C$, that is, there exists a number $\eta > 0$ such that 
\be
\label{etaineq}
\sup_{C \in {\cal C}} \lambda(\partial(C : \pi) ) 
 \ > \eta \, \mbox{ for every finite measurable partition } \, \pi . 
\ee
Using the inequality (\ref{etaineq}) we construct a sequence of ``splitting sets'' 
$S_1, S_2, \ldots \subseteq \X$ from the sets in $\C$ in a stage-wise
fashion.  
At the $k$th stage the splitting set $S_k$ is obtained from 
a sequential procedure that makes use of the splitting sets 
$S_1,\ldots, S_{k-1}$ produced at previous stages.  The splitting
sets are used to identify arbitrarily large finite collections of sets in 
$\C$ having full join.   The existence of these collections implies
that $\C$ has infinite VC dimension by Lemma \ref{Join}. 

\vskip.1in

\noindent
{\bf First stage.}
Define the refining sequence of joins $J_1(n) = C_1 \vee \cdots \vee C_n$ 
for $n \geq 1$.  It follows from (\ref{etaineq}) that for each $n$ there is a set
$C_1(n) \in \C$ whose boundary $G_1(n) = \partial(C_1(n) : J_1(n))$ has 
measure greater than $\eta$.  
Note that the sets $\{ G_1(n) : n \geq 1 \}$ are measurable $\sigS_0$.  
By standard results in functional analysis, there exists 
a subsequence $\{ n_m \}$ and an $\sigS_0$-measurable function $h_1$
such that $\int g \, I_{G_1(n_m)} \, d\mu \to \int g \, h_1 \, d\mu$ as
$m$ tends to infinity for every $g \in L_2(\X, \sigS_0,\mu)$.  (The 
function $h_1$ is the weak limit of the indicator functions $I_{G_1(n_m)}$.)
It follows that
$0 \leq h_1 \leq 1$ almost surely, and that $\int h_1 \, d\mu \geq \eta$.
Define the splitting set $S_1 = \{ h_1 > 0 \}$ and note that $\mu(S_1) \geq \eta$.

For simplicity, let $J_1(m)$, $C_1(m)$, and $G_1(m)$ denote, respectively, the quantities
$J_1(n_m)$, $C_1(n_m)$, and $G_1(n_m)$ along the subsequence defining $h_1$.  We
adopt similar notation for subsequences encountered at subsequent stages.

\vskip.2in

\noindent
{\bf Subsequent stages.}
Suppose now that we have constructed splitting sets $S_j$ at stages
$j = 1,\ldots, k-1$, and wish
to construct the splitting set $S_k$ at stage $k$.  
Begin by defining the refining sequence of joins 
$J_k(n) \ = \ S_1 \vee \cdots \vee S_{k-1} \vee C_1 \vee \cdots \vee C_n$
for $n \geq 1$.  It follows from (\ref{etaineq}) that for each $n$ 
there is a set $C_{k}(n) \in \C$ whose boundary $G_{k}(n) = \partial(C_k(n) : J_k(n))$ 
has measure greater than $\eta$.  Proceeding as in Stage 1, there is a subsequence 
$\{ I_{G_k(m)} \}$ having a weak limit $h_k \in L_2(\X, \sigS_0, \mu)$ such that
$0 \leq h_k \leq 1$ almost surely, and $\int h_k \, d\mu \geq \eta$.
Define the splitting set $S_k = \{ h_k > 0 \}$ and note that $\mu(S_k) \geq \eta$.

\vskip.2in

\noindent
{\bf Construction of Full Joins.}
Fix an integer $L \geq 1$.   As the measure of each splitting set $S_k$ 
is at least $\eta$, there exist positive integers 
$k_1 < k_2 < \ldots < k_{L+1}$ 
such that $\mu(\bigcap_{j=1}^{L+1} S_{k_j}) > 0$.   
Suppose for simplicity, and without loss of generality, that $k_j = j$. 
For $l = 1, \ldots, L+1$ define 
\[
Q_l = \bigcap_{j=1}^{l} S_j
\] 
In what follows we will make repeated use of the elementary fact that 
$\int_B (h_1 \cdots h_l) \, d\mu > 0$
if and only if $\mu(B \cap Q_l) > 0$.

We claim that there exist sets $D_1,\ldots,D_L \in \C$ such that for
each $l = 1, \ldots, L$,
\be
\label{positive}
\int_B (h_1 \cdots h_l) \, d\mu > 0
\ \mbox{ for every } \ 
B \in D_l \vee \cdots \vee D_L .
\ee
The inequalities (\ref{positive}) are established by reverse induction, 
beginning with the case $l = L$.  To this end, note that 
\[
0
\  < \  
\int (h_1 \cdots h_{L+1}) \, d\mu
\ = \ 
\lim_{m \to \infty} \int (h_1 \cdots h_L) I_{G_{L+1}(m)} \, d\mu ,
\]
and therefore $\mu( Q_L \cap G_{L+1}(m) ) > 0$ for all $m$ sufficiently
large. Fix such an $m$ and let $D = C_{L+1}(m)$.  It follows from the definition of 
$G_{L+1}(m)$ that for some cell 
$A \in J_{L+1}(m)$,
\be
\label{relns1}
\mu( Q_L \cap A ) > 0 \ \mbox{ and } \ 
\mu(A \cap D) \cdot \mu(A \cap D^c) > 0 .
\ee
The inclusion of the sets $S_1,\ldots, S_L$ in the definition of the
joins $J_{L+1}(n)$ ensures that $Q_L$ is a finite union of cells
of $J_{L+1}(m)$.  
The first relation in (\ref{relns1}) then implies that $A$ is necessarily 
a subset of $Q_L$, and it follows from the second 
relation that $\mu(Q_L \cap D) \cdot \mu(Q_L \cap D^c) > 0$.  Letting
$D_L = D$ the last inequality implies (\ref{positive}) in the case
$l = L$.


Suppose now that for some $1 < l < L$ we have identified sets 
$D_l, D_{l+1}, \ldots, D_L$ such that (\ref{positive}) holds.  
Then for each cell $B$ in the join $D_l \vee \cdots \vee D_L$,
\[
0
\  < \  
\int_B (h_1 \cdots h_l) \, d\mu
\ = \ 
\lim_{m \to \infty} \int_B (h_1 \cdots h_{l-1}) I_{G_l(m)} \, d\mu .
\]
Therefore, there exists an integer $m$ such that
$\mu(B \cap Q_{l-1} \cap G_l(m)) > 0$ for every
$B \in D_l \vee \cdots \vee D_L$.
As the join $J_l(m)$ includes the first $n_m$ elements of $\C$, by
enlarging $m$ if necessary we may assume that $J_l(m)$ includes
$D_l, \ldots, D_L$.  Let $D = C_l(m)$ and 
let $B$ be any cell of $D_l \vee \cdots \vee D_L$.  
The definition of $G_l(m)$ implies that for some cell 
$A \in J_l(m)$,
\be
\label{relns2}
\mu(B \cap Q_{l-1} \cap A) > 0 
\ \mbox{ and } \
\mu(A \cap D) \cdot \mu(A \cap D^c) > 0 .
\ee
Both $Q_{l-1}$ and $B$ are
equal to a union of cells of the partition $J_1(m)$, so
the first relation in (\ref{relns2}) implies that $A \subseteq B \cap Q_{l-1}$,
and it then follows from the second relation that
$\mu(B \cap Q_{l-1} \cap D)$ and $\mu(B \cap Q_{l-1} \cap D^c)$ 
are positive.
As these inequalities hold for each $B \in D_l \vee \cdots \vee D_L$, we have
$\int_{B'} (h_1 \cdots h_{l-1}) \, d\mu > 0$
for every $B' \in D \vee D_l \vee \cdots \vee D_L$.  Letting $D_{l-1} = D$
completes the induction.

It follows from (\ref{positive}) that the sets $D_1,\ldots,D_L$ have full join,
and as $L \geq 1$ was arbitrary, Lemma \ref{Join} implies that $\C$ has
infinite VC dimension, which completes the proof of the theorem.

\vskip.2in

\noindent
{\bf Remark:}
An inspection of the proof shows that the approximating partitions $\pi$ in the 
theorem can be taken to be measurable $\sigma(\C)$.
A simple counterexample shows that $\pi$ may not be chosen from the smaller family 
$\bigcup_{n=1}^{\infty} \sigma(C_1 \vee C_2\vee \ldots \vee C_n)$. 
Let $\X=[0,1]$ and let $\mu$ be Lebesgue measure. Let $a_1, a_2, \ldots$ 
be a sequence of positive real numbers such that $s=\sum_{n=1}^{\infty} a_n < 1$. 
Define $s_0=0$ and $s_n = \sum_{i=1}^n a_i$ for $n \geq 1$, and let 
$C_n = [s_{n-1},s_n)$. 
Clearly, the VC-dimension of the class $\{C_1,C_2,\ldots\}$ 
equals 1, since the sets are disjoint.
Define $J_n = C_1\vee C_2\vee \ldots \vee C_n$. 
Then the set $A_n=[s_n,1]$ is a
single element in $J_n$ with measure 
$1-s_n > 1-s > 0$.  Moreover, both $A_n\cap C_{n+1}$ 
and $A_n\cap C^{\prime}_{n+1}$ have positive measure.  
Thus, for $n \geq 1$, $A_n \subseteq \partial(C_{n+1}:G_n)$
and $\mu (\partial(C_{n+1}:G_n)) > 1-s$.

\vskip.3in

\noindent
{\bf\Large  Acknowledgements} \\
The authors are indebted to an anonymous referee of the earlier paper
\cite{AdNob10a} who suggested the general form of Theorem \ref{BdyThm2}, 
and whose detailed comments led to a simpler and more 
general proof.  
The authors would also like to acknowledge helpful discussions with 
Ramon van Handel, who provided
feedback on an earlier version of this work \cite{AdNob10c}, and who brought
the papers \cite{BourFremTala78, Rosen74, Tala84} to our attention.
The work presented in this paper was supported in part
by NSF grant DMS-0907177.


\begin{thebibliography}{9}
\csname bibmessage\endcsname

\bibitem{AdNob10a}
\textsc{Adams, T.M.} and \textsc{Nobel, A.B.} (2010)
Uniform convergence of Vapnik-Chervonenkis classes under 
ergodic sampling.
\textit{Annals of Probability} 
\textbf{38}(4)1345-1367.

\bibitem{AdNob10b}
\textsc{Adams, T.M.} and \textsc{Nobel, A.B.} (2010)
The gap dimension and uniform laws of large numbers for ergodic processes. 
Preprint. arXiv:1007.2964v1

\bibitem{AdNob10c}
\textsc{Adams, T.M.} and \textsc{Nobel, A.B.} (2010)
Uniform approximation and bracketing properties of VC classes.
Preprint. arXiv1007.4037v1
		
\bibitem{Assouad83}
\textsc{Assouad, P.} (1983)
Densit\'e et dimension.
\textit{Annales de l'Institut Fourier} 
\textbf{33}(3) 233-282.
MR0723955 (86j:05022)	
	
	





\bibitem{BourFremTala78}
\textsc{Bourgain, J.} and \textsc{Fremlin, D.H.} and \textsc{Talagrand, M.} (1978) 
Pointwise compact sets of Baire measurable functions.
\textit{American Journal of Mathematics} 
\textbf{100} 845-886.
MR0509077 (80b:54017)




\bibitem{DevGyoLug96}
\textsc{Devroye, L.} and \textsc{Gy\"orfi, L.} and \textsc{Lugosi, G.} (1996) 
\textit{A probabilistic theory of pattern recognition}.
Springer. MR1383093 (97d:68196)

\bibitem{Dudley99}
\textsc{Dudley, R.M.} (1999) 
\textit{Uniform Central Limit Theorems.}
Cambridge University Press, Cambridge.
MR1720712 (2000k:60050)

\bibitem{GaenStute76}
\textsc{Gaenssler, P.} and \textsc{Stute, W.} (1976)
On uniform convergence of measures with applications to uniform convergence of empirical distributions.
\textit{Empirical distributions and processes (Selected Papers, Meeting on 
Math. Stochastics, Oberwolfach, 1976)} 
45-56.  Springer Lecture Notes in Math., \textbf{566}.
MR0433534  (55 \#6510)





\bibitem{KeaSch94}
\textsc{Kearns, M.J.} and \textsc{Schapire, R.E.} (1994)
Efficient distribution-free learning of probabilistic concepts.
\textit{Journal of Computer and System Sciences}
\textbf{48(3)} 464--497. 


\bibitem{Mat02}
\textsc{Matousek, J.} (2002)
Lectures on Discrete Geometry.
\textit{Graduate Texts in Mathematics} 
\textbf{212} Springer, New York. MR1899299 (2003f:52011)





%

%

\bibitem{Poll84}
\textsc{Pollard, D.} (1984)
\textit{Convergence of Stochastic Processes}
Springer, New York. MR0762984 (86i:60074)


\bibitem{Rosen74}
\textsc{Rosenthal, H.P.} (1974) 
A characterization of Banach spaces containing l1
\textit{Proceedings of the National Academy of Sciences U.S.A.} 
\textbf{71} 2411-2413.
MR0358307 (50 \#10773)


%

\bibitem{Tala84}
\textsc{Talagrand, M.} (1984)
\textit{Pettis integral and measure theory} 
Memoirs of the American Mathematics Society  
\textbf{51}(307). MR0756174  (86j:46042)

\bibitem{Tala87}
\textsc{Talagrand, M.} (1987)
The Glivenko-Cantelli problem. 
\textit{Annals of Probability} 
\textbf{15:3} 837--870. MR0893902 (88h:60012)


\bibitem{VaaWel96}
\textsc{van der Vaart, A.W.} and \textsc{Wellner, J.A.} (1996)
\textit{Weak Convergence and Empirical Processes}.
Springer-Verlag, New York. 
MR1385671 (97g:60035)

\bibitem{Vapnik00}
\textsc{Vapnik, V.N.} (2000)
\textit{The nature of statistical learning theory}.  Second edition.
Springer-Verlag, New York. 
MR1719582 (2001c:68110)

\bibitem{VC71}
\textsc{Vapnik, V.N.} and \textsc{Chervonenkis, A.Ya.} (1971)
On the uniform convergence of relative frequencies of events to their probabilities.
\textit{Theory of Probability and its Applications}
\textbf{16} 264--280. MR0627861 (83d:60031)


\end{thebibliography}
\end{document}